\documentclass[reqno]{amsart}

\usepackage{amsmath}
\usepackage{amssymb}
\usepackage{amsfonts}
\usepackage{mathrsfs}
\usepackage{bbm}
\usepackage{graphicx}

\usepackage{todonotes}

\graphicspath{ {./Pictures/} }

\usepackage[tmargin=1.5in,bmargin=1.5in,lmargin=1.3in,rmargin=1.3in]{geometry}
\usepackage{microtype}

\newtheorem{thm}{Theorem}

\theoremstyle{remark}

\newcommand\nlongleftrightarrow{\mathrel{\,\,\,\not\!\!\!\longleftrightarrow}}

\title[Sharpness of the Phase Transition for the Corrupted Compass Model]{Sharpness of the Phase Transition for the Corrupted Compass Model on Transitive Graphs}

\author{Thomas Beekenkamp}
\address{Mathematisches Institut, Ludwig-Maximilians-Universit\"at M\"unchen, Theresienstra\ss{}e 39, 80333 M\"unchen, Germany}
\email{Thomas.Beekenkamp@math.lmu.de}
\date{\today}

\begin{document}
\begin{abstract}
In the corrupted compass model on a vertex-transitive graph, a neighbouring edge of every vertex is chosen uniformly at random and opened. Additionally, with probability $p$, independently for every vertex, every neighbouring edge is opened. We study the size of open clusters in this model. Hirsch et al. \cite{hirsch2018} have shown that for small $p$ all open clusters are finite almost surely, while for large $p$, depending on the underlying graph, there exists an infinite open cluster almost surely. We show that the corresponding phase transition is sharp, i.e., in the subcritical regime, all open clusters are exponentially small. Furthermore we prove a mean-field lower bound in the supercritical regime. The proof uses the by now well established method using the OSSS inequality. A second goal of this note is to showcase this method in an uncomplicated setting.

\end{abstract}
\maketitle

\section{Introduction and Main Result}
\begin{figure}\label{fig:cc}
  \centering
    \includegraphics[width=0.95\textwidth]{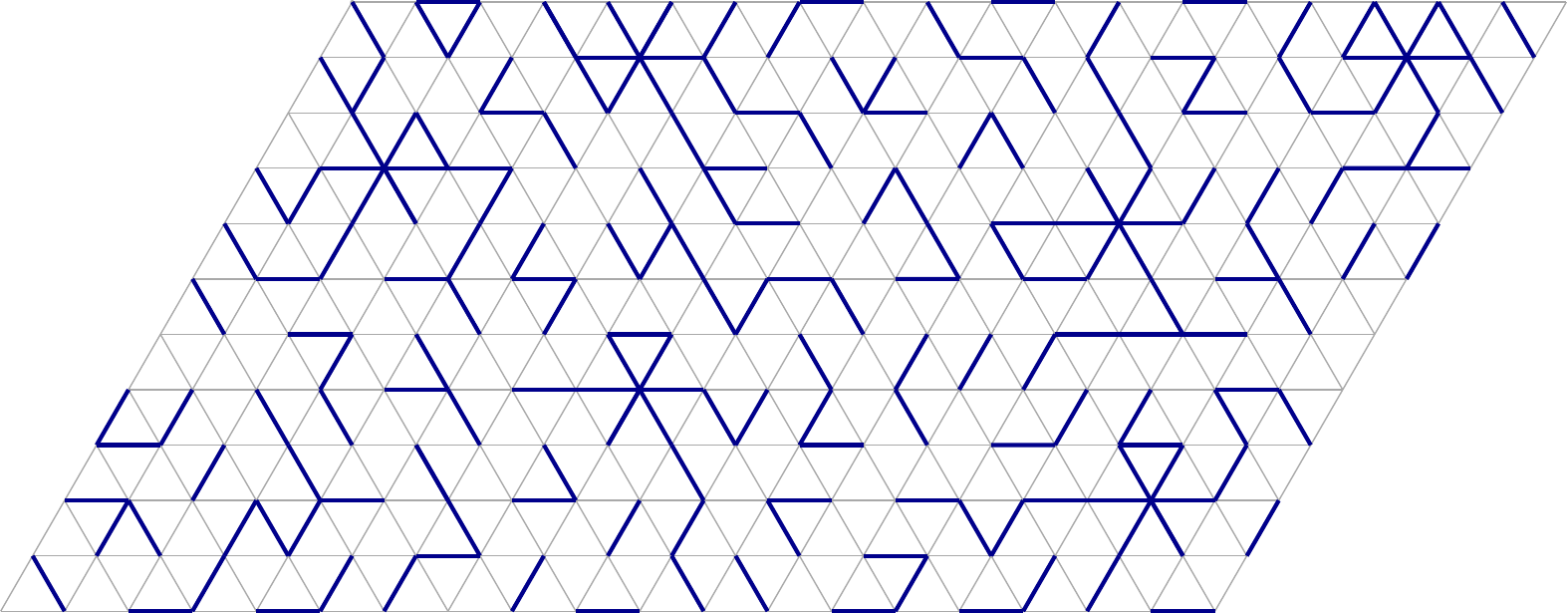}
    \caption{The corrupted compass model on the triangular lattice.}  
\end{figure}

Let $G=(V,E)$ be an infinite, connected, locally finite, vertex-transitive graph. We consider the corrupted compass model on $G$, which is informally defined as follows. Each vertex $v\in V$ is corrupted with probability $p$, independently of each other. For each corrupted vertex, we declare each neighbouring edge to be open. On the other hand, for an uncorrupted vertex, we choose one neighbouring edge to be open uniformly at random. A possible configuration of this model on the triangular lattice is shown in Figure \ref{fig:cc}.

The corrupted compass model was introduced by Hirsch, Holmes and Kleptsyn \cite{hirsch2018} in the context of reinforcement models for neural networks. They show that in a class of reinforcement models the reinforced edges almost surely do not form an infinite cluster if the reinforcement is strong enough. They show this by making a coupling between the reinforcement model and the corrupted compass model, and subsequently showing that in the latter model there exists only finite clusters almost surely for $p$ small enough.

The corrupted compass model is not only relevant to reinforcement models, as the model was also used in the context of alignment percolation by Beaton, Grimmett and Holmes \cite{beaton2019}. In the one-choice alignment percolation model on $\mathbb{Z}^d$ introduced by these authors, a Bernoulli site percolation configuration with parameter $p$ is taken. Subsequently, for each occupied vertex, one of the $2d$ directions is chosen uniformly at random the entire line segment in this direction until the next occupied vertex is declared blue. The authors then ask the question whether there exists an infinite blue cluster. The main problem in the analysis of this model is the lack of monotonicity in $p$. Nevertheless, the authors show that for $p$ large enough there exists no infinite blue clusters almost surely. They show this by dominating the alignment percolation model by a corrupted compass model with parameter $1-p$. Since the corrupted compass model does not have infinite clusters for $1-p$ small enough, the one-choice alignment percolation model does not have any infinite clusters for $p$ large enough.

Our contribution is to show that the phase transition in the corrupted compass model is sharp. That is, in the subcritical regime the open clusters are exponentially small. Sharp phase transitions are a common theme in percolation models. The sharpness of the phase transition for Bernouilli percolation was first shown by Menshikov \cite{menshikov1986} and independently by Aizenman and Barsky \cite{AizenmanBarsky1987}. While being different, both proofs have proven to be hard to generalise to models with dependencies. In particular, it is not clear how to apply these proofs to the corrupted compass model, since the status of neighbouring edges is dependent. A technique to prove sharp phase transitions was introduced by Duminil-Copin, Raoufi and Tassion \cite{Duminil-Copin2019} using the OSSS inequality for Boolean functions. This technique, which we will also use, has proven to be fruitful for the random cluster model \cite{Duminil-Copin2019}, Voronoi percolation \cite{Duminil-Copin2017}, Boolean percolation \cite{Duminil-Copin2018}, the Widom-Rowlinson model \cite{Dereudre2018}, and level sets of Gaussian fields \cite{Duminil-Copin2020}\cite{Muirhead2020}.

To precisely state our results, we first introduce several definitions. For $v\in V$, let $\mathcal{N}(v)$ denote the set of edges that include $v$, and let $d=|\mathcal{N}(v)|$ (which is independent of $v$). We fix an arbitrary vertex $0\in V$ to be the origin. For $v,w\in V$, let $d(v,w)$ denote the graph distance between $v$ and $w$ in $G$. For $n\in \mathbb{N}$, we define the balls
\[
\Lambda_n^v:=\{w\in V \::\: d(v,w)\leq n\}, \quad\quad \partial \Lambda_n^v:=\{w\in V \::\: d(v,w)= n\}.
\] 
For $v=0$, we drop part of the notation: $\Lambda_n=\Lambda_n^0$ and $\partial \Lambda_n=\partial \Lambda_n^0$. For a bond configuration $\eta$ and $v,w\in V$, we say that $v\longleftrightarrow w$, if there is a path of open edges starting in $v$ and ending in $w$. Similarly, for $A\subset V$ we say that $v\longleftrightarrow A$, whenever there exists $w\in A$, such that $v\longleftrightarrow w$. We say that $0\longleftrightarrow\infty$, if for all $n\in \mathbb{N}$, we have $0\longleftrightarrow\partial \Lambda_n$. We define the critical value for percolation as
\[
p_c:=\sup\{p\::\: \mathbb{P}_p(0\longleftrightarrow\infty)=0\}.
\]

Hirsch, Holmes and Kleptsyn \cite{hirsch2018} have shown that, for $p$ small enough, all clusters are finite almost surely. From this it follows that $p_c>0$. On the other hand the corrupted compass model dominates the Bernoulli site percolation model that only uses the corrupted compasses. Therefore, we have $p_c\leq p_c^{\text{site}}(G)$, where $p_c^{\text{site}}(G)$ is the critical threshold for Bernoulli site percolation on $G$. Depending on the graph $G$, this threshold is nontrivial, so that also $p_c<1$. 

\newpage
In this note we will prove the following theorem.
\begin{thm}\label{thm}
Consider the corrupted compass model with parameter $p$.
\begin{enumerate}
\item Let $p<p_c$. There exists a constant $c>0$, such that for all $n\in \mathbb{N}$,
\[
\mathbb{P}_p(0\longleftrightarrow\partial \Lambda_n)\leq \exp(-cn).
\]
\item There exists a constant $c>0$, such that for all $p>p_c$,
\[
\mathbb{P}_{p}(0\longleftrightarrow\infty)\geq c(p-p_c).
\]
\end{enumerate}

\end{thm}

\section{The OSSS inequality and sharp phase transitions}
We now precisely define the corrupted compass model on $G$. We consider the probability space $(\Omega, \mathcal{F},\mathbb{P}_p)$, where
\[
\Omega= \prod_{v\in V} \,[0,1]\times \mathcal N(v),
\]
the $\sigma$-algebra $\mathcal{F}$ is generated by the cylindrical events, and $\mathbb{P}_p$ is the product measure of the uniform measures on $[0,1]\times \mathcal N(v)$. For $\omega\in \Omega$, we denote $U_v:=\omega_{v,1}$, i.e., the uniform random variable on $[0,1]$ associated to $v$, and $A_v:=\omega_{v,2}$, the uniformly chosen edge in  $\mathcal{N}(v)$. We define $X_v:=(U_v,A_v)$. Let $\mathcal{K}$ denote the set of corrupted vertices, i.e., 
\[
\mathcal{K}:= \{v\in V\::\: U_v<p\}.
\]
We can obtain the bond configuration $\eta$ as follows. Let $\eta:\Omega \to \{0,1\}^E$ be given by 
\[
\eta_e(\omega):=\mathbbm{1}\left\{e\in \bigcup_{v\in \mathcal K}\mathcal{N}(v) \,\cup\, \bigcup_{v\in \mathcal{K}^c} \{A_v\} \right\}.
\]
We say that an edge $e$ is open whenever $\eta_e=1$, and closed otherwise.

We now introduce the framework for the OSSS inequality for Boolean functions. For simplicity, and since it is sufficient for our ends, we only consider Boolean functions on finite domains. For $n\in \mathbb{N}$, let
\[
\Omega_n= \bigotimes_{v\in \Lambda_n} \,[0,1]\times \mathcal N(v).
\]
A function $f:\Omega_n \to \{0,1\}$ is called a Boolean function. Since $(\Omega_n,\mathcal{F},\mathbb{P}_p)$ is a probability space, the function $f$ is also a random variable, and we can talk about the variance of $f$. The OSSS inequality is a bound on the variance of a Boolean function. In order to introduce the inequality, we need the notion of a decision tree. Let $T$ be a decision tree that determines the value of $f$. That is, $T$ starts by revealing $X_{v_0}$ for some (possibly random) $v_0\in \Lambda_n$. Depending on the value of $X_{v_0}$, it chooses $v_1\in \Lambda_n$ and reveals the value of $X_{v_1}$. Then, depending on the information obtained so far, it again chooses another vertex $v_2\in \Lambda_n$ and reveals $X_{v_2}$. This process goes on until $T$ has enough information to determine the value of $f$, i.e., until the values of the unrevealed variables cannot change the value of $f$ anymore. We define the revealment of $X_v$ to be
\[
\text{Rev}_v(T):=\mathbb{P}_p( T \text{ reveals } X_v).
\]
The influence of $X_v$ on $f$ is defined as
\[
\text{Inf}_v(f):=\mathbb{P}_p\big(\mathbbm{1}\{f(\omega)\neq f(\tilde{\omega}_v)\big),
\]
where $\tilde{\omega}_v$ is obtained from $\omega$ by resampling $X_v$ independently. The OSSS inequality states that for any decision tree $T$, we have
\[
\text{var}(f)\leq \sum_{v} \text{Rev}_v(T)\text{Inf}_v(f).
\]
This inequality was proven by O'Donnell , Saks, Schramm and Servedio \cite{ODonnell2005}. It is an improvement on the Poincar\'e inequality, which states that $\text{var}(f)\leq \sum_v \text{Inf}_v(f)$, and is obtained from basic Fourier analysis for Boolean functions. The OSSS discounts the influences of the vertices that are unlikely to be seen by the decision tree. For a broader introduction to analysis of Boolean functions, see O'Donnell \cite{ODonnell2014}.

The general approach to proving sharpness of a phase transition in random spatial models is to find a suitable differential inequality. In particular, differential inequalities can be found for the one arm event, i.e., for the quantity $\theta_n(p):=\mathbb{P}_p(0\longleftrightarrow\partial \Lambda_n)$. The original proof by Menshikov for Bernoulli percolation uses a differential inequality for $\theta_n(p)$, but this proof relies heavily on the independencies of the model. Alternatively, differential inequalities can be found for the magnetisation, a particular quantity derived from the cluster size distribution. In the proof of Aizenman and Barsky for Bernoulli percolation, the authors find two differential inequalities for the magnetisation, also by harvesting the independencies of the model. For dependent models a more robust approach is needed. 

The OSSS inequality can be used to prove similar differential inequalities in dependent settings. One way to deal with the dependencies is by writing the probability as a product space, as we did for $\Omega_n$. This is possible if the model exhibits a sufficient amount of independence, such as in our case. The proofs for Boolean percolation and Voronoi percolation also take this approach. Another possible way to deal with dependencies is possible when the model admits a monotonic measure. In particular this is the case for the random cluster model, a generalization of Bernoulli percolation and the Ising model. Hutchcroft proved a differential inequality for the magnetisation in the random cluster model using the OSSS inequality, from which the sharpness of the phase transition follows. The original proof by Duminil-Copin, Raoufi and Tassion for the sharpness in the random cluster model also uses the fact that the measure is monotonic, but they obtain a differential inequality for the one arm event.

We will also apply the OSSS inequality to the one arm event, that is, to the Boolean function $f:=\mathbbm{1}\{0\longleftrightarrow\partial \Lambda_n\}$, to find a differential inequality for $\theta_n(p)$. This is in rough lines the same approach as in the original proof by Duminil-Copin, Raoufi and Tassion. However, the precise way to apply the OSSS inequality varies from model to model, since a suitable decision tree has to be chosen. As the name suggests, the decision tree should be able to decide whether $0$ is connected to $\partial \Lambda_n$. Moreover, it should do so by looking at as little variables as possible, so that on average the revealment is small and hence the differential inequality is strong. One option is to explore from $\partial \Lambda_k$, for $1\leq k\leq n$ chosen uniformly at random. In this way, the average revealment is small, and it is determined if $0\longleftrightarrow \partial \Lambda_n,$ since such a connection has to go through $\partial \Lambda_k$. The sum of the influences can then typically be bounded by $\theta'_n(p)$ using a type of Russo's formula, so that we indeed obtain a differential inequality. The precise way to do this is also model dependent. We will now carry out this procedure for the corrupted compass model. 

\section{Proof}
\begin{figure}
  \centering
    \includegraphics[width=0.7\textwidth]{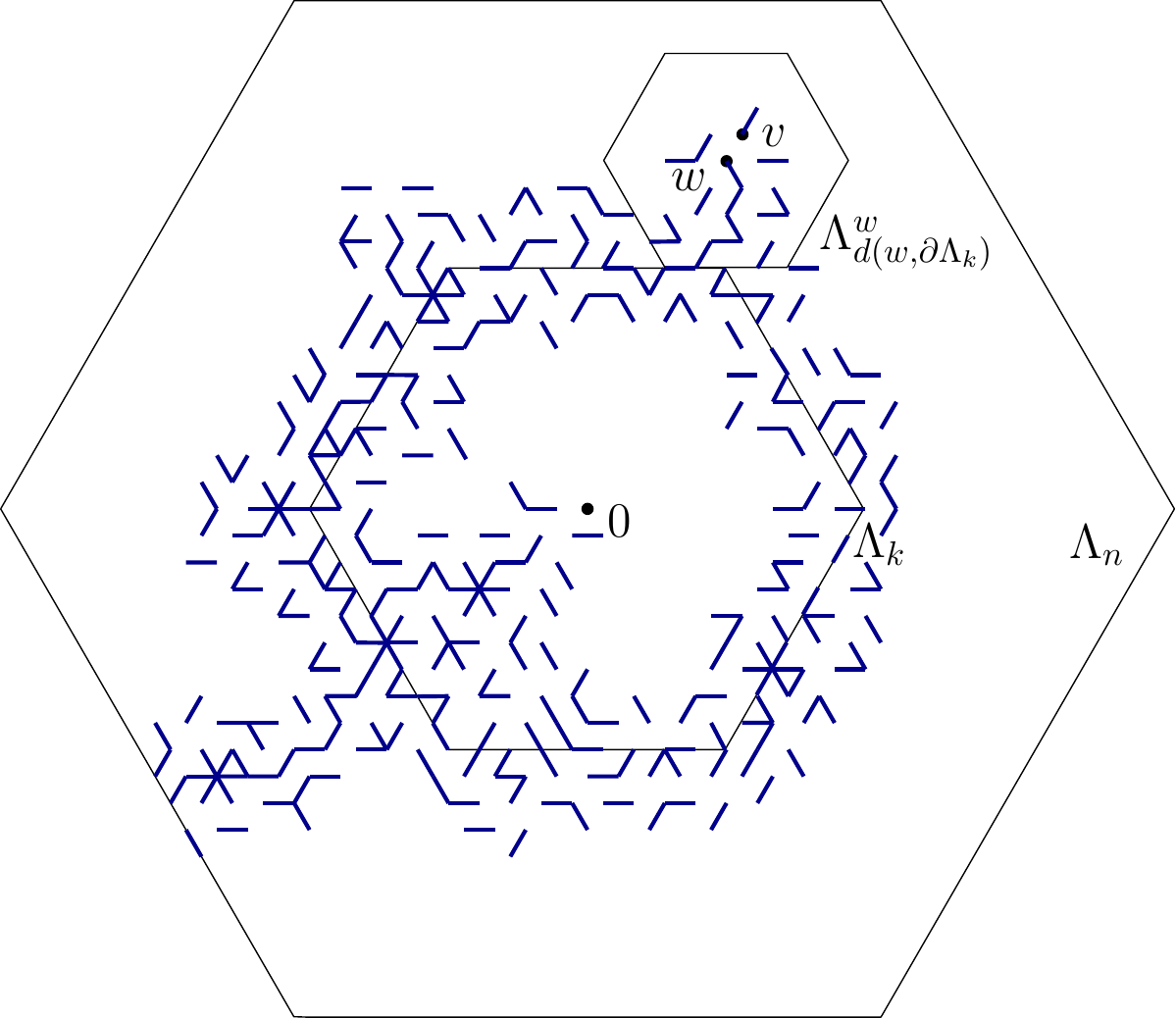}
    \caption{The decision tree $T_k$ exploring the cluster of $\partial \Lambda_k$. When $X_v$ is revealed, there must be a neighbour $w$ of $v$ that is connected to $\partial \Lambda_k$.} 
    \label{fig:cc2}
\end{figure}

We will apply the OSSS inequality to the Boolean function $f:=\mathbbm{1}\{0\longleftrightarrow\partial \Lambda_n\}$. This function only depends on the variables in  $\{X_v : v\in \Lambda_n\}$. We write $\theta_n(p):=\mathbb{P}_p(0\longleftrightarrow\partial \Lambda_n)$. For $1\leq k\leq n$, let $T_k$ be the decision tree that explores the cluster of $\partial \Lambda_k$. This decision tree determines $f$, since a path from $0$ to $\partial \Lambda_n$ must go through $\partial \Lambda_k$. To precisely describe the decision tree $T_k$, we need a subalgorithm $\textsf{Determine}(v)$, for $v\in \Lambda_n$. When $\textsf{Determine}(v)$ is called, $X_v$ is revealed, as well as $X_w$ for all neighbours $w$ of $v$. This determines the state of the edges in $\mathcal{N}(v)$. The decision tree $T_k$ now does the following. Firstly, $\textsf{Determine}(v)$ is called for all $v\in \partial \Lambda_k$. This determines the vertices in $\partial \Lambda_{k-1}$ and $\partial \Lambda_{k+1}$ that are connected to $\partial \Lambda_k$ by open edges. For these vertices again, $\textsf{Determine}(v)$ is called to see which vertices are connected to them. This process continues until the entire cluster of $\partial \Lambda_k$ inside $\Lambda_n$ is determined. In particular, it is then determined whether there is a connection from $0$ to $\partial \Lambda_n$. The exploration process carried out by $T_k$ for the model on the triangular lattice is shown in Figure \ref{fig:cc2}

Applying the OSSS inequality to $f$ and $T_K$, and summing over $k$ gives
\begin{equation}
n\theta_n(p)(1-\theta_n(p))\leq \sum_{v\in \Lambda_n}\sum_{k=1}^n\text{Rev}_v(T_k)\text{Inf}_v(f).
\end{equation}

\subsection{Bound on the Revealment}
By summing over $k$, we essentially average over all spheres $\partial \Lambda_k$ with radius up to $n$, so that the average revealment is small. This is in spirit the same as taking $1\leq k\leq n$ uniformly at random. We note that if $X_v$ is revealed by $T_k$, it follows that $\Lambda_1^v\longleftrightarrow \partial \Lambda_k$. We obtain
\begin{align*}
\sum_{k=1}^n \text{Rev}_v(T_k)\leq \sum_{k=1}^n \mathbb{P}_p(\Lambda_1^v\longleftrightarrow \partial \Lambda_k)&\leq \sum_{k=1}^n \sum_{w\in \Lambda_1^v} \mathbb{P}_p(w\longleftrightarrow \partial \Lambda_k)\\
&\leq \sum_{k=1}^n \sum_{w\in \Lambda_1^v} \mathbb{P}_p\big(w\longleftrightarrow \partial \Lambda^w_{d(w,\partial \Lambda_k)}\big).
\end{align*}
Using translation invariance, we have
\[
\sum_{k=1}^n \mathbb{P}_p\big(w\longleftrightarrow \partial \Lambda^w_{d(w,\partial \Lambda_k)}\big)\leq 2\sum_{k=1}^n \mathbb{P}_p(0\longleftrightarrow\partial \Lambda_k).
\]
If we define $S_n=S_n(p):=\sum_{k=1}^n\theta_k(p)$, it follows that
\[
\sum_{k=1}^n \text{Rev}_v(T_k)\leq 2d S_n.
\]

\subsection{Bound on the Influence}
For $\omega\in \Omega$, we say that $v\in V$ is a pivotal corrupted compass for an event $A$, whenever $\mathbbm{1}_A(\omega)\neq \mathbbm{1}_A(\hat{\omega}_v)$, where $\hat{\omega}_v$ is obtained from $\omega$ by corrupting $v$ if $v$ is uncorrupted in $\omega$, or by uncorrupting $v$ if $v$ is corrupted in $\omega$. Russo's formula gives
\[
\theta_n'(p)=\frac{\mathrm{d}}{\mathrm{d}p}\mathbb{P}_p(0\longleftrightarrow\partial\Lambda_n)=\sum_{v\in \Lambda_n} \mathbb{P}_p(v \text{ pivotal corrupted compass for }0\longleftrightarrow\partial \Lambda_n).
\]
The aim is to relate the above quantity to the total influence, so that we obtain a differential inequality. We have
\[
\sum_{v\in V} \text{Inf}_v(f)= 2\sum_{v\in V}\mathbb{P}_p(f(\omega)=0, f(\tilde{\omega}_v)=1),
\]
where $\tilde{\omega}_v$ is obtained from $\omega$ by resampling $X_v$ independently. If $f(\omega)=0$, but $f(\tilde{\omega}_v)=1$, it follows that $v$ is not corrupted in $\omega$. Therefore, corrupting $v$ will put $f$ to $1$, because this will open at least as much edges as the resampling of $X_v$. Thus $v$ is a pivotal corrupted compass. We obtain
\[
\sum_{v\in V} \text{Inf}_v(f)\leq 2\sum_{v\in V}\mathbb{P}_p(v \text{ pivotal corrupted compass for }0\longleftrightarrow\partial \Lambda_n).
\]
Hence,
\[
\theta_n'(p)\geq \frac{1}{2}\sum_{v\in V}\text{Inf}_v(f).
\]

Combining the OSSS inequality and the bounds on the revealment and the influence gives
\begin{equation}
\frac{\mathrm{d}}{\mathrm{d}p}\theta_n(p)\geq \frac{n}{4dS_n}\theta_n(p)(1-\theta_n(p)).
\end{equation}

\subsection{Proof of Theorem \ref{thm}}
To finish the proof, we distinguish between the cases $p_c=1$ and $p_c<1$. First we assume that $p_c=1$. Let $p_0<p_c$. We have $\theta_n(p_0)\to 0$ as $n\to \infty$. Let $N$ be such that $\theta_n(p_0)\leq\tfrac{1}{2}$ for all $n>N$. Then for all $p \leq p_0$ and for all $n>N$ we have
\[
\theta_n'(p)\geq \frac{n}{4dS_n}\theta_n(p)(1-\theta_n(p))\geq \frac{n}{4dS_n}\theta_n(p)(1-\theta_n(p_0))\geq \frac{1}{8d}\frac{n}{S_n}\theta_n(p).
\]
From this inequality we can obtain the sharpness of the phase transition, which we will show in the next section. First we will find the same differential inequality, but with a different constant for the case $p_c<1$. We can assume that $d\geq 3$, since the only infinite, connected, transitive graph with $d=2$ is $\mathbb{Z}$ with nearest neighbour edges, for which $p_c=1$. Let $p_c<\delta<1$. For $n\geq 2$ and $p\leq \delta$, we bound
\[
1-\theta_n(p)\geq 1-\theta_2(\delta)=\mathbb{P}_\delta(0\nlongleftrightarrow\partial \Lambda_2).
\]
We can construct a configuration in which $0\nlongleftrightarrow\partial \Lambda_2$, and which has positive probability, as follows. Let $v$ be the vertex that the compass of $0$ points to, i.e., $A_0=\{0,v\}$. We require that the compass of $v$ points back to $0$, which happens with probability $1/d$. Furthermore we want $0$, $v$, and all other neighbours of $0$ and $v$ to be uncorrupted, which costs at most $(1-\delta)^{2d}$. Finally, we want that the compasses of the other neighbours of $0$ and $v$ do not point towards $0$ or $v$, which happens with probability $(\tfrac{d-2}{d})^{2d-2}$. All together we find
\[
\mathbb{P}_\delta(0\nlongleftrightarrow\partial \Lambda_2)\geq (1-\delta)^{2d}\frac{1}{d}\left(\frac{d-2}{d}\right)^{2d-2}=:C_0>0,
\]
so that for all $n\geq 2$ and all $p\leq \delta$, we have
\begin{equation}\label{eq:difineq}
\theta_n'(p)\geq C_1\frac{n}{S_n}\theta_n(p),
\end{equation}
where $C_1:=C_0/4d$. Since $C_1\leq 1/8d$, the above inequality holds for the case where $p_c=1$ as well, for $n>N$.

\subsection{Analysis of the differential inequality}
The remainder of the proof of Theorem \ref{thm} involves analysing the above differential inequality (\ref{eq:difineq}). Suppose $p<p_c$. We will first show the exponential decay of $\theta_n(p)$. For that purpose, let $p<p_1<p_2<p_c$. From (\ref{eq:difineq}), it follows that
\[
\log \theta_n(p) '\geq C_1\frac{n}{S_n},
\]
for all $n>N$. Integrating the above inequality from $p_1$ to $p_2$ gives
\[
-\log \theta_n(p_1)\geq \log \theta_n(p_2)-\log \theta_n(p_1)\geq C_1(p_2-p_1)\frac{n}{S_n(p_2)},
\]
so that
\begin{equation}\label{eq:difana}
\theta_n(p_1)\leq \exp\left(-C_1(p_2-p_1)\frac{n}{S_n(p_2)}\right).
\end{equation}

If $S_n(p_2)$ is bounded in $n$, i.e., if $\sum_{k=0}^\infty\theta_k(p_2)$ converges, the desired exponential decay would follow from the above inequality. In fact, it suffices if $S_n(p_2)\leq n^{1-\alpha}$ for $0<\alpha<1$ and $n$ large enough: from (\ref{eq:difana}) it then follows that 
\[
\theta_n(p_1)\leq \exp\left(-C_1(p_2-p_1)n^\alpha\right),
\]
for $n$ large enough, so that $\sum_{k=0}^\infty\theta_k(p_1)$ converges. We can then bootstrap this result by using the inequality (\ref{eq:difana}) again to find the desired exponential decay. This motivates the definition of the following critical point, which we will show to be equal to $p_c$:
\[
\tilde{p}_c:=\sup\{p\::\: \limsup_{n\to \infty} \frac{\log S_n}{\log n}<1\}.
\]
If $p_2<\tilde{p}_c$, there exists $0<\alpha<1$ such that $S_n(p_2)\leq n^{1-\alpha}$ for $n$ large enough. It then follows that we have stretched exponential decay at $p_1$ and exponential decay at $p$.

It remains to show that there exists $c>0$ such that linear lower bound $\theta(p)\geq c(p-\tilde{p}_c)$ holds for all $p>\tilde{p}_c$, since from this the equality $p_c=\tilde{p}_c$ follows. Let $\tilde{p}_c<p_1<p$. Using (\ref{eq:difineq}), we find
\begin{align*}
\sum_{k=1}^n\frac{\theta_k(p)'}{k}\geq C_1\sum_{k=1}^n \frac{\theta_k(p)}{S_k}\geq C_1 \sum_{k=1}^n\int_{S_k}^{S_{k+1}}\frac{1}{t}\,\mathrm{d}t&=C_1 \sum_{k=1}^n (\log S_{k+1}-\log S_k)\\
&=C_1(\log S_{n+1}-\log S_1)\geq C_1\log S_{n+1}
\end{align*}
We define $T_n(p):=\frac{1}{\log n}\sum_{k=1}^n \frac{\theta_k(p)}{k}$, and find
\[
T_n(p)'\geq C_1\frac{\log S_{n+1}}{\log n}.
\]
Integrating the above inequality from $p_1$ to $p$ gives
\[
T_n(p)-T_n(p_1)\geq C_1(p-p_1)\frac{\log S_{n+1}(p_1)}{\log n}.
\]
Note that for all $p$, $\frac{1}{\log n}T_n(p)\to \theta(p)$ for $n\to \infty$, so that
\[
\theta(p)\geq\theta(p)-\theta(p_1)=\limsup_{n\to \infty} (T_n(p)-T_n(p_1))\geq C_1 (p-p_1)\limsup_{n\to \infty}\frac{\log S_{n+1}(p_1)}{\log n}\geq C_1 (p-p_1)>0,
\]
since $p_1>\tilde{p}_c$. Because $p>\tilde{p}_c$ is arbitrary, it follows that $\tilde{p}_c=p_c$. By letting $p_1\to p_c$, we find the desired lower bound.

\qed

\bibliographystyle{abbrv}
\bibliography{cc}

\end{document}